\documentclass[12pt]{article}
\usepackage{latexsym,amsmath,amssymb}
\input epsf

\usepackage[utf8]{inputenc}

\usepackage{amssymb}

\renewcommand{\ge}{\geqslant}
\renewcommand{\le}{\leqslant}

\title{\textbf{Power sum expansion of chromatic quasisymmetric functions}}

\author{Christos~A.~Athanasiadis\\
\small Department of Mathematics\\[-0.8ex]
\small University of Athens\\[-0.8ex]
\small Athens 15784, Hellas (Greece)\\[-0.8ex]
\small\tt caath@math.uoa.gr
}

  \def\QQ{{\mathbb Q}}
  \def\ZZ{{\mathbb Z}}

  \def\aA{{\mathcal A}}
  \def\cC{{\mathcal C}}
  \def\dD{{\mathcal D}}
  
  \def\nN{{\mathcal N}}
  \def\oO{{\mathcal O}}

  \def\asc{{\rm asc}}
  \def\Des{{\rm Des}}

  \def\QSym{{\rm QSym}}
  
  \def\inv{{\rm inv}}
  
  \def\SYT{{\rm SYT}}
   
  \def\sm{\smallsetminus}
  \newcommand{\qed}{$\hfill \Box$}

\date{\small April 3, 2015}

\begin{document}
\maketitle

\newtheorem{theorem}{Theorem}[section]
\newtheorem{proposition}[theorem]{Proposition}
\newtheorem{corollary}[theorem]{Corollary}
\newtheorem{defn}[theorem]{Definition}
\newtheorem{remark}[theorem]{Remark}
\newtheorem{lemma}[theorem]{Lemma}
\newtheorem{example}[theorem]{Example}
\newtheorem{examples}[theorem]{Examples}
\newtheorem{conjecture}[theorem]{Conjecture}
\newtheorem{fact}[theorem]{Fact}
\newtheorem{question}[theorem]{Question}
\newtheorem{observation}[theorem]{Observation}
\newtheorem{claim}[theorem]{Claim}

\begin{abstract}
The chromatic quasisymmetric function of a graph was introduced by
Shareshian and Wachs as a refinement of Stanley's chromatic symmetric function. An explicit combinatorial formula, conjectured by Shareshian and Wachs, expressing the chromatic quasisymmetric function of the incomparability graph of a natural unit interval order in terms of power sum symmetric functions, is proven. The proof uses a formula of Roichman for the irreducible characters of the symmetric group.

\end{abstract}

\section{Introduction}
\label{sec:intro}

The chromatic quasisymmetric function of a graph was introduced by
Shareshian and Wachs \cite{SW12, SW14} as a refinement of
Stanley's symmetric function generalization of the chromatic
polynomial of a graph \cite{Sta95}. To recall the definition,
consider a simple graph $G$ on a totally ordered set of vertices
$V = \{v_1, v_2,\dots,v_n\}$. We denote by $\cC(G)$ the set of
proper colorings of $G$, meaning the set of maps $\kappa: V \to
\{1, 2,\dots\}$ such that $\kappa(v_i) = \kappa(v_j)$ implies
$\{v_i, v_j\}$ is not an edge of $G$. For such a coloring
$\kappa$, we denote by $\asc(\kappa)$ the number of edges $\{v_i,
v_j\}$ of $G$ with $i < j$ for which $\kappa(v_i) < \kappa (v_j)$.
Given a sequence $x = (x_1, x_2,\dots)$ of commuting
indeterminates, the chromatic quasisymmetric function of $G$ is
defined as

  \begin{equation} \label{eq:defchqf}
    X_G (x, t) \ = \ \sum_{\kappa \in \cC(G)} t^{\asc(\kappa)} \,
                     x_\kappa
  \end{equation}
where $x_\kappa = x_{\kappa(v_1)} x_{\kappa(v_2)} \cdots x_{\kappa(v_n)}$.

\smallskip
The function $X_G$ is quasisymmetric in $x$; its study \cite{SW12, SW14} connects seemingly disparate topics, such as graph colorings, permutation statistics and the cohomology of Hessenberg varieties,
and provides valuable insight into Stanley's chromatic symmetric function \cite{Sta95}, to which $X_G$ reduces for $t=1$. When $G$
is the incomparability graph of a natural unit interval order, the function $X_G$ turns out to be symmetric in $x$ \cite[Theorem~4.5]{SW14}. Moreover, its setting allows for interesting refinements
of several theorems and conjectures on the expansion of Stanley's chromatic symmetric function in various bases of the algebras of symmetric and quasisymmetric functions. The present note confirms
one of the conjectures of \cite{SW12, SW14} (see
\cite[Conjecture~4.15]{SW12} \cite[Conjecture~7.6]{SW14}),
predicting an explicit combinatorial formula which expresses $X_G$
in terms of power sum symmetric functions. This formula refines a formula of Stanley \cite[Theorem~2.6]{Sta95} and provides some evidence in favor of the positivity conjecture of Shareshian and 
Wachs (see \cite[Conjecture~4.9]{SW12} \cite[Conjecture~1.3]{SW14})
on the expansion of $X_G$ in the basis of elementary symmetric functions.

Our main result (Theorem~\ref{thm:main}) is formulated and proven in Section~\ref{sec:main}, after background material has been reviewed in Section~\ref{sec:pre}. The result is derived from the expansion of 
$X_G$ in the basis of fundamental quasisymmetric functions \cite[Theorem~3.1]{SW14} via a formula of Roichman~\cite{Roi97} for the irreducible characters of the symmetric group (the proof is motivated by the discussion of character formulas and distributions of descent sets in \cite[Section~6]{AR13}). Section~\ref{sec:orient} exploits a known bijection \cite{Stei92} between certain permutations of a partially ordered set $P$ and acyclic orientations of the incomparability graph of $P$ to provide an alternative formulation (Corollary~\ref{cor:main}) of the main result.

\section{Preliminaries}
\label{sec:pre}

This section fixes notation and reviews the background material from \cite{AR13, SW12, SW14} which is necessary to state and prove our results. Any unexplained notation and terminology on partial orders and symmetric and quasisymmetric functions can be found in \cite[Chapter~3]{StaEC1} \cite[Chapter~7]{StaEC2}. We denote by $\mathfrak{S}(S)$ the group of permutations, and by $|S|$ the cardinality, of a finite set $S$. For nonnegative integers $n$ we set $[n] := \{1, 2,\dots,n\}$ and write $\mathfrak{S}_n$ instead of $\mathfrak{S}([n])$.

\subsection{Symmetric and quasisymmetric functions}
\label{subsec:sym}

Given a sequence $x = (x_1, x_2,\dots)$ of commuting
indeterminates and a commutative ring (with unit) $R$, we will
denote by $\Lambda^n_R$ (respectively, by $\QSym^n_R$) the
$R$-module of homogeneous symmetric (respectively, quasisymmetric)
functions of degree $n$ over $R$ in $x$. The \emph{fundamental
quasisymmetric function} associated to $S \subseteq [n-1]$ is
defined as
  $$ F_{n, S} \ = \ \sum_{\substack{ i_1 \le i_2 \le \cdots \le i_n \\ j \in S \Rightarrow i_j < i_{j+1} }} x_{i_1} x_{i_2} \cdots x_{i_n} $$
(this definition is slightly different from the one used in
\cite{SW12, SW14}). The set $\{ F_{n, S}: S \subseteq [n-1]\}$ is known to be a basis of $\QSym^n_\ZZ$.

There are two involutions on $\QSym^n_R$, denoted by $\omega$ and
$\rho$, which will be of interest here. They can be defined by letting $\omega(F_{n, S}) = F_{n, [n-1] \sm S}$ and $\rho(F_{n, S}) = F_{n, n-S}$ for $S \subseteq [n-1]$, where $n - S := \{ n - a: a \in S\}$, and extending by linearity. The involution $\omega$ restricts to the standard involution on $\Lambda^n_R$; see \cite[Section~7.6]{StaEC2}. The involution $\rho$, discussed in \cite[Section~2]{SW14}, fixes every element of $\Lambda^n_R$.

We now recall two ways to express a Schur function in terms of fundamental quasisymmetric functions and power sum symmetric functions. We denote by $\SYT (\lambda)$ the set of standard Young tableaux of shape $\lambda \vdash n$. For such a tableau $Q$, we denote by $\Des(Q)$ the set of $i \in [n-1]$ for which $i+1$ appears in $Q$ in a lower row than $i$. Stanley's theory of $P$-partitions implies (see \cite[Theorem~7.19.7]{StaEC2}) that
 \begin{equation} \label{eq:sFexpansion}
    s_\lambda(x) \ = \ \sum_{Q \in \SYT (\lambda)} F_{n, \Des(Q)}
  \end{equation}
where, as usual, $s_\lambda(x)$ is the Schur function associated to $\lambda$.

We denote by $\chi^\lambda$ the irreducible character of the symmetric group $\mathfrak{S}_n$ associated to the partition $\lambda \vdash n$ and by $p_\lambda (x)$ the corresponding power sum symmetric function. We also set $z_\lambda = \prod_{i \ge 1} i^{m_i} \, m_i!$ where $m_i$ is the number of parts of $\lambda$ which are equal to $i$.

\smallskip
For any sequence $\alpha = (\alpha_1, \alpha_2,\dots,\alpha_k)$ of positive integers summing to $n$ (called composition of $n$), we write $S(\alpha) = \{r_1, r_2,\dots,r_k\}$ for the set of partial sums $r_i = \alpha_1 + \alpha_2 + \cdots + \alpha_i$ of $\alpha$ and set $r_0 = 0$. Following \cite[Definition~3.1]{AR13}, we call \emph{$\alpha$-unimodal} any subset of $[n-1]$ whose intersection with each of the sets $\{r_{i-1} + 1,\dots,r_i - 1\}$ for $1 \le i \le k$ is a prefix (possibly empty) of the latter. For instance, $S(\alpha)$ is $\alpha$-unimodal since all these intersections are empty. Moreover, if $\alpha = (n)$, then the $\alpha$-unimodal subsets of $[n-1]$ are those equal to $[p]$ for some $p \in \{0, 1,\dots,n-1\}$ and if $\alpha = (1, 1,\dots,1)$, then every subset of $[n-1]$ is $\alpha$-unimodal. As another example, if $\alpha = (3, 1, 4, 2)$, then $\{1, 3, 5, 6\}$ is $\alpha$-unimodal but $\{1, 3, 5, 7\}$ is not.

We denote by $U_\alpha$ the set of $\alpha$-unimodal subsets of $[n-1]$. The following statement is a special case of \cite[Theorem~4]{Roi97}.

\begin{theorem} {\rm (\cite[Theorem~3.6]{AR13})}  \label{thm:roi}
For all partitions $\lambda, \mu \vdash n$
  \begin{equation} \label{eq:chiroi}
    \chi^\lambda (\mu) \ = \ \sum_{\substack{Q \in \SYT(\lambda), \\ 
                                   \Des(Q) \in U_\mu }}
                        (-1)^{|\Des(Q) \sm S(\mu)|}.
  \end{equation}
As a result,
  \begin{equation} \label{eq:schroi}
    s_\lambda (x) \ = \ \sum_{\mu \vdash n} z^{-1}_\mu p_\mu (x)
                        \sum_{\substack{ Q \in \SYT (\lambda), \\ 
                              \Des(Q) \in U_\mu }}
                        (-1)^{|\Des(Q) \sm S(\mu)|}
  \end{equation}
in $\Lambda^n_\QQ$ for every $\lambda \vdash n$.
\end{theorem}

\subsection{Chromatic quasisymmetric functions}
\label{subsec:chrom}

We will need two results from \cite{SW12, SW14} on the chromatic function $X_G$, so we first recall some definitions. Let $G = (V, E)$ be a simple graph on the vertex set $V = \{v_1, v_2,\dots,v_n\}$, totally ordered by $v_1 < v_2 < \cdots < v_n$. Given a sequence $\sigma = (\sigma_1, \sigma_2,\dots,\sigma_p)$ of pairwise distinct elements of $V$, we denote by $\inv_G(\sigma)$ the number of \emph{$G$-inversions} of $\sigma$, meaning pairs $(i, j)$ of indices with $1 \le i < j \le p$, such that $\sigma_i > \sigma_j$ and $\{\sigma_i, \sigma_j\} \in E$.
Assume that $G$ is the incomparability graph of a partial order $P$ on the set $V$, meaning the simple graph with vertices the elements of $V$ and edges the pairs of elements which are incomparable in $P$. We say that an index $i \in [n-1]$ is a \emph{$P$-descent} of $\sigma$ if $\sigma_i >_P \sigma_{i+1}$ and denote by $\Des_P (\sigma)$ the set of $P$-descents of $\sigma$. We say that $\sigma_j$ is a
\emph{left-to-right $P$-maximum} of $\sigma$ if $\sigma_i <_P \sigma_j$ for all indices $i < j$.

\begin{theorem} {\rm (\cite[Theorem~3.1]{SW14})}  \label{thm:SWexpF}
Let $G$ be the incomparability graph of a poset $P$ on a totally ordered set of cardinality $n$. Then
  \begin{equation} \label{eq:SWexpF}
    \omega X_G (x, t) \ = \ \sum_{\sigma \in \mathfrak{S} (P)}
                   t^{\inv_G(\sigma)} \, F_{n, n - \Des_P (\sigma)}.
  \end{equation}
\end{theorem}

A partial order $P$ on the set $[n]$ is said to be a \emph{natural unit interval order} if there exist real numbers $y_1 < y_2 < \cdots < y_n$ such that $i <_P j \Leftrightarrow y_i + 1 < y_j$ for all $i, j \in [n]$; several equivalent characterizations are given in \cite[Section~4]{SW14}. Clearly, the class of natural unit interval orders is closed under taking induced subposets. Shareshian and Wachs showed \cite[Theorem~4.5]{SW14} that if $G$ is the incomparability graph of a natural unit interval order on the set $[n]$, endowed with its natural total order, then $X_G$ is a symmetric function in $x$.

We will denote by $\nN(P)$ the set of permutations $\sigma \in \mathfrak{S}(P)$ which (written in one line notation) have no
$P$-descent and no nontrivial (meaning, other than the initial term) left-to-right $P$-maximum.

\begin{lemma} {\rm (\cite[Lemma~7.4]{SW14})}  \label{lem:SWbn}
Let $G$ be the incomparability graph of a natural unit interval order $P$ on the set $[n]$ and let $\omega X_G (x, t) = \sum_{\lambda \vdash n} z^{-1}_\lambda b_\lambda(t) p_\lambda (x)$ be the expansion of $\omega X_G$ in the basis of power sum symmetric functions. Then
  \begin{equation} \label{eq:SWbn}
    b_{(n)}(t) \ = \ \sum_{\sigma \in \nN(P)} t^{\inv_G(\sigma)}.
  \end{equation}
\end{lemma}

Lemma~\ref{lem:SWbn} is shown in~\cite{SW14} using the 
Murnaghan-Nakayama rule and a formula for the expansion of 
$X_G$~\cite[Theorem~6.3]{SW14} in the basis of Schur functions. A more elementary proof of this lemma appears in 
Section~\ref{sec:orient} (see Remark~\ref{rem:alt}).

\section{Power sum expansion}
\label{sec:main}

To state our main result, we need to introduce one more piece of notation from \cite{SW12, SW14}. Given a partially ordered set $P$ with $n$ elements and a partition $\lambda = (\lambda_1, \lambda_2,\dots,\lambda_k)$ of $n$, we denote by $\nN_\lambda (P)$ the set of all permutations $\sigma \in \mathfrak{S}(P)$ such that when we break
$\sigma$ (written in one line notation) from left to right into contiguous segments of lengths $\lambda_1, \lambda_2,\dots,\lambda_k$, each contiguous segment has neither a $P$-descent, nor a nontrivial left-to-right $P$-maximum.

The following theorem appeared as \cite[Conjecture~4.15]{SW12} and as part of \cite[Conjecture~7.6]{SW14} (the unimodality part of
\cite[Conjecture~7.6]{SW14} remains open). It has been verified in a special case in \cite[Proposition~7.9]{SW14}; for an application, see \cite[Proposition~9.12]{SW14}.

\begin{theorem} \label{thm:main}
Let $P$ be a natural unit interval order on the set $[n]$ with incomparability graph $G$. Then
  \begin{equation} \label{eq:pexp}
    \omega X_G (x, t) \ = \ \sum_{\lambda \vdash n} z^{-1}_\lambda
                          p_\lambda (x)
                          \sum_{\sigma \in \nN_\lambda(P)}
                          t^{\inv_G(\sigma)}.
  \end{equation}
\end{theorem}

The proof will be deduced from the following statement, which is perhaps implicit in the discussion in \cite[Section~9]{AR13} (see Proposition 9.3 there).

\begin{proposition} \label{prop:AR}
Let $X(x) \in \Lambda^n_R$ be a homogeneous symmetric function of degree $n$ over a commutative $\QQ$-algebra $R$ and suppose that
  \begin{equation} \label{eq:propARa}
    X(x) \ = \ \sum_{S \subseteq [n-1]} a_S F_{n,S}
  \end{equation}
for some $a_S \in R$. Then
  \begin{equation} \label{eq:propARb}
    X(x) \ = \ \sum_{\lambda \vdash n} z^{-1}_\lambda p_\lambda (x)
            \sum_{S \in U_\lambda} (-1)^{|S \sm S(\lambda)|} \, a_S,
  \end{equation}
where $U_\lambda$ is the set of $\lambda$-unimodal subsets of
$[n-1]$.
\end{proposition}

\noindent
\emph{Proof.} 
We may write $X(x) = \sum_{\lambda \vdash n} c_\lambda s_\lambda
(x)$ for some coefficients $c_\lambda \in R$. Replacing $s_\lambda
(x)$ with the right-hand side of (\ref{eq:schroi}) for $\lambda 
\vdash n$ we get $X(x) = \sum_{\mu \vdash n} b_\mu z^{-1}_\mu p_\mu 
(x)$, where
  \begin{eqnarray}
    b_\mu &=& \sum_{\lambda \vdash n} c_\lambda
                \sum_{\substack{ Q \in \SYT(\lambda), \\ 
                      \Des(Q) \in U_\mu }}
                (-1)^{|\Des(Q) \sm S(\mu)|} \label{eq:ARdefb}
      \nonumber \\
  & & \nonumber \\
  &=& \sum_{S \in U_\mu} (-1)^{|S \sm S(\mu)|} \, \sum_{\lambda   
      \vdash n} \, c_\lambda \cdot |\{Q \in \SYT (\lambda): \Des(Q) 
      = S\}|. \label{eq:ARb}
  \end{eqnarray}
Using Equation~(\ref{eq:sFexpansion}) to expand $X(x)$ in the basis of fundamental quasisymmetric functions gives
  $$ X(x) \ = \ \sum_{\lambda \vdash n} c_\lambda s_\lambda (x) \ = \
             \sum_{\lambda \vdash n} c_\lambda \sum_{Q \in \SYT
             (\lambda)} F_{n, \Des(Q)}. $$
Comparing the previous equation with (\ref{eq:propARa}) we get
 \begin{equation} \label{eq:ARac}
    a_S \ = \ \sum_{\lambda \vdash n} c_\lambda \cdot |\{Q \in \SYT
               (\lambda): \Des(Q) = S\}|.
  \end{equation}
From (\ref{eq:ARb}) and (\ref{eq:ARac}) we conclude that $b_\mu = \sum_{S \in U_\mu} (-1)^{|S \sm S(\mu)|} \, a_S$ for every $\mu \vdash n$ and the proof follows.
\qed

\bigskip
\noindent
\emph{Proof of Theorem~\ref{thm:main}}.
We recall that the involution $\rho$ on $\QSym^n_R$, discussed in
Section~\ref{subsec:sym}, fixes each element of $\Lambda^n_R$. Since $\omega X_G (x, t) \in \Lambda^n_{\QQ [t]}$, we may deduce from
(\ref{eq:SWexpF}) that
  $$ \omega X_G (x, t) \ = \ \rho (\omega X_G (x, t)) \ = \
     \sum_{\sigma \in \mathfrak{S}_n} t^{\inv_G(\sigma)} \, \rho (
     F_{n, n - \Des_P (\sigma)}) \ = \ \sum_{\sigma \in \mathfrak{S}_n}
     t^{\inv_G(\sigma)} \, F_{n, \Des_P (\sigma)}. $$
This expression and Proposition~\ref{prop:AR} imply that

  \begin{equation} \label{eq:Xbp}
    \omega X_G (x, t) \ = \ \sum_{\lambda \vdash n} z^{-1}_\lambda
                             b_\lambda(G, t) \, p_\lambda (x)
  \end{equation}
where

  \begin{equation} \label{eq:b_lambda}
    b_\lambda(G, t) \ = \ \sum_{\substack{ \sigma \in \mathfrak{S}_n
                          \\ \Des_P(\sigma) \in U_\lambda }}
    (-1)^{|\Des_P(\sigma) \sm S(\lambda)|} \, t^{\inv_G (\sigma)}
  \end{equation}
for every $\lambda \vdash n$.

We fix a partition $\lambda =
(\lambda_1,\dots,\lambda_k)$ of $n$ and note that there is a
one-to-one correspondence from $\mathfrak{S}_n$ to the set of tuples $(\pi, \sigma_1,\dots,\sigma_k)$, where $\pi = (B_1,\dots,B_k)$ is an ordered partition of $[n]$ with block sizes $|B_i| = \lambda_i$ (say that such an ordered partition has type $\lambda$) and $\sigma_i \in \mathfrak{S}(B_i)$ for every index $i$. Specifically, the $\sigma_i$ are the segments of $\sigma \in \mathfrak{S}_n$ when the latter (written in one line notation) is broken from left to right into contiguous segments of sizes $\lambda_1,\dots,\lambda_k$. Clearly, we have $\Des_P(\sigma) \in U_\lambda$ if and only if $\Des_P(\sigma_i) \in U_{(\lambda_i)}$ for every index $i$. Moreover, if $\sigma \in \mathfrak{S}_n$ corresponds to $(\pi, \sigma_1,\dots,\sigma_k)$, then
  \begin{equation} \label{eq:inv}
    \inv_G (\sigma) \ = \ \inv_G(\pi) \ + \ \sum_{i=1}^k \inv_G
                          (\sigma_i),
  \end{equation}
where $\inv_G(\pi)$ stands for the number of $G$-inversions of $\sigma$ which involve two elements in different blocks of $\pi$, and $\Des_P(\sigma) \sm S(\lambda)$ is equal to the disjoint union of the sets $\Des_P(\sigma_i)$ for $1 \le i \le k$. The previous remarks imply that (\ref{eq:b_lambda}) may be rewritten as
  \begin{equation} \label{eq:b_lambda2}
    b_\lambda(G, t) \ = \ \sum_{\pi = (B_1,\dots,B_k) \in
    \Pi_\lambda} t^{\inv_G (\pi)} \ \prod_{i=1}^k
    \sum_{\substack{ \sigma_i \in \mathfrak{S}(B_i)
    \\ \Des_P(\sigma_i) \in U_{(\lambda_i)} }}
    (-1)^{|\Des_P(\sigma_i)|} \, t^{\inv_G (\sigma_i)},
  \end{equation}
where $\Pi_\lambda$ is the set of all ordered partitions of $[n]$ of type $\lambda$.

The special case $\lambda = (n)$ of Equation~(\ref{eq:b_lambda}), combined with Lemma~\ref{lem:SWbn}, implies that

 \begin{equation} \label{eq:b(n)}
    \sum_{\substack{ \sigma  \in \mathfrak{S}_n
    \\ \Des_P(\sigma) \in U_{(n)} }}
    (-1)^{|\Des_P(\sigma)|} \, t^{\inv_G (\sigma)} \ = \
    \sum_{\sigma \in \nN(P)} t^{\inv_G(\sigma)}.
  \end{equation}
We recall that the induced subposet of $P$ on each block $B_i$ is (under the natural order-isomorphism of the latter with the set $[\lambda_i]$) again a natural unit interval order. Thus, applying (\ref{eq:b(n)}) to these subposets, we may rewrite (\ref{eq:b_lambda2}) as
  \begin{equation} \label{eq:b_lambda3}
    b_\lambda(G, t) \ = \ \sum_{\pi = (B_1,\dots,B_k) \in
    \Pi_\lambda} t^{\inv_G (\pi)} \ \prod_{i=1}^k \, \left(
    \sum_{\sigma_i \in \nN(B_i)} t^{\inv_G(\sigma_i)} \right).
  \end{equation}
Using the inverse of the correspondence mentioned after
(\ref{eq:b_lambda}) and, in view of (\ref{eq:inv}), we find that Equation~(\ref{eq:b_lambda3}) yields the desired expression
$b_\lambda(G, t) = \sum_{\sigma \in \nN_\lambda(P)} t^{\inv_G (\sigma)}$ for the coefficient of $z^{-1}_\lambda p_\lambda (x)$ in
$\omega X_G$.
\qed

\section{Permutations and acyclic orientations}
\label{sec:orient}

This section exploits a bijection due to Steingr\'imsson~\cite{Stei92} and provides an alternative formulation of Theorem~\ref{thm:main} in terms of acyclic orientations of the incomparability graph $G$. As another application, a short proof Lemma~\ref{lem:SWbn} is given.

We begin with a few definitions. An \emph{acyclic orientation} of a simple graph $G = (V, E)$ is a choice of direction $a \to b$ or $b \to a$ for each edge $\{a, b\} \in E$ so that the resulting directed graph has no directed cycle. A \emph{sink} of such an orientation is a vertex $a \in V$ such that no edge of $G$ is directed away from $a$. Given an acyclic orientation $o$ of $G$, we may define a partial order $\bar{o}$ on the vertex set $V$ by requiring that for $a, b \in V$ we have $a \le_{\bar{o}} b$ if there exists a directed path in
$G$ with initial vertex $b$ and terminal vertex $a$ (note that the minimal elements of $\bar{o}$ are exactly the sinks of $o$). We denote by $\aA\oO(G)$ the set of acyclic orientations of $G$.

Assume that $G$ is the incomparability graph of a partially ordered set $P$ with $n$ elements and denote by $\dD(P)$ the set of permutations $\sigma \in \mathfrak{S}(P)$ which have no $P$-descent. Every permutation $\sigma \in \mathfrak{S}(P)$ gives rise to an acyclic orientation $\varphi_P (\sigma)$ for which an edge $\{a, b\}$ of $G$ is directed as $b \to a$ if $a$ precedes $b$ in $\sigma$ (written in one line notation). Steingr\'imsson showed \cite[Theorem~4.12]{Stei92} (see also \cite[Exercise~3.60]{StaEC1} \cite[Lemma~4.6]{Sta02}) that the restriction $\varphi_P: \dD(P) \to \aA\oO(G)$ of the resulting map to $\dD(P)$ is bijective for every finite poset 
$P$. The inverse map can be described as follows. Given an acylic orientation $o \in \aA\oO(G)$, note that the sinks of $o$ are pairwise nonadjacent in $G$ and hence they form a chain in $P$. We denote by $\sigma_1$ the smallest element of this chain, define the elements $\sigma_2,\dots,\sigma_n$ by the same procedure from the poset obtained by deleting $\sigma_1$ from $P$ and set $\psi_P (o) = (\sigma_1, \sigma_2,\dots,\sigma_n)$.
For $i \in [n-1]$ the element $\sigma_i$ is the $P$-minimum sink of the restriction of $o$ on the induced subgraph of $G$ on the vertex set $\{\sigma_i, \sigma_{i+1},\dots,\sigma_n\}$ and similarly for $\sigma_{i+1}$. This implies that either $\sigma_i <_P \sigma_{i+1}$,  or else $\{ \sigma_i, \sigma_{i+1} \}$ is an edge of $G$ directed as $\sigma_{i+1} \to \sigma_i$ by $o$; in particular, we cannot have $\sigma_i >_P \sigma_{i+1}$. We conclude that $\psi_P (o)$ has no 
$P$-descents and hence we have a well defined map $\psi_P: \aA \oO(G) \to \dD(P)$.

The following theorem summarizes the properties of the maps $\varphi_P$ and $\psi_P$ which are important for us. The last statement assumes that the vertex set of $G$ has been equipped with an arbitrary total order. For $o \in \aA\oO(G)$, we then denote by $\asc(o)$ the number directed edges $b \to a$ of $G$ for which $a$ is larger than $b$ in this total order.

\begin{lemma} \label{lem:Stei}
The maps $\varphi_P: \dD(P) \to \aA\oO(G)$ and $\psi_P: \aA \oO(G) \to \dD(P)$ are bijections which are inverses of each other. Moreover, for every $\sigma \in \dD(P)$ the left-to-right $P$-maxima of $\sigma$ are exactly the sinks of $\varphi_P (\sigma)$ and $\inv_G(\sigma) = \asc(\varphi_P(\sigma))$.
\end{lemma}

\noindent
\emph{Proof.} 
A proof that $\psi_P (\varphi_P(\sigma)) = \sigma$ for every $\sigma \in \dD(P)$ can be extracted from the proof of \cite[Theorem~4.12]{Stei92}; see, for instance, the proof of \cite[Lemma~4.6]{Sta02}. That $\varphi_P(\psi_P(o)) = o$ for every $o \in \aA\oO(G)$ follows from the fact that $\psi(o)$ is, by construction, a linear extension of the poset $\bar{o}$. The claim that $\inv_G(\sigma) = \asc(\varphi_P(\sigma))$ for every $\sigma \in \dD(P)$ follows from the definition of $\varphi_P$. Finally, given $\sigma \in \dD(P)$, the same definition shows that an element $b \in P$ is a sink of $\varphi_P (\sigma)$ if and only if $b$ is comparable in $P$ to every element preceding $b$ in $\sigma$. Since $\sigma$ has no $P$-descent, the latter can happen only when $b$ is a left-to-right $P$-maximum of $\sigma$.
\qed

\bigskip
We recall from Section~\ref{sec:main} that $\Pi_\lambda$ denotes the set of ordered partitions of $[n]$ of type $\lambda$. For a graph $G = ([n], E)$ and $\pi = (B_1, B_2,\dots,B_k) \in \Pi_\lambda$ we denote by $\aA\oO^\ast (G, \pi)$ the set of (necessarily acyclic) orientations $o$ of $G$ which satisfy the following: (a) $o$ restricts to an acyclic orientation with a unique sink on each induced subgraph of $G$ with vertex set $B_i$; and (b) every edge $\{a, b\} \in E$ with $a \in B_i$ and $b \in B_j$ for some indices $i < j$ is directed by $o$ as $b \to a$. The following statement is a consequence of Equation (\ref{eq:b_lambda3}) and
Lemma~\ref{lem:Stei}.

\begin{corollary} \label{cor:main}
Let $P$ be a natural unit interval order on the set $[n]$ with incomparability graph $G$. Then $\omega X_G (x, t) = \sum_{\lambda \vdash n} z^{-1}_\lambda b_\lambda(G, t) \, p_\lambda (x)$ where

  \begin{equation} \label{eq:pexpalt}
    b_\lambda (G, t) \ = \ \sum_{\pi \in \Pi_\lambda} \,
                           \sum_{o \in \aA\oO^\ast (G, \pi)}
                           t^{\asc(o)}
  \end{equation}
for every $\lambda \vdash n$.
\end{corollary}

\begin{remark} \label{rem:alt}
\rm As mentioned in Section~\ref{subsec:chrom}, the proof of 
Lemma~\ref{lem:SWbn} given in~\cite{SW14} uses the Murnaghan-Nakayama rule and a formula for the expansion of $X_G$~\cite[Theorem~6.3]{SW14} in the basis of Schur functions. A proof which does not assume these results can be given as follows. As noted in \cite[Section~7]{SW14}, the coefficient $b_{(n)} (t)$, considered in this lemma, is equal to the coefficient of $e_n (x)$ when $X_G (x, t)$ is expanded in the basis of elementary symmetric functions. It follows from the proof of \cite[Theorem~3.3]{Sta95} (see \cite[Theorem~5.3]{SW14}) that this coefficient is equal to the sum of $t^{\asc(o)}$, where $o$ ranges over all acyclic orientations of $G$ with a unique sink. Equation~(\ref{eq:SWbn}) follows from this statement and 
Lemma~\ref{lem:Stei}.
\end{remark}

\end{document}